\documentclass[english]{article}
\usepackage[T1]{fontenc}
\usepackage[latin9]{inputenc}
\usepackage{amsmath}
\usepackage{amssymb}
\usepackage{tikz}

\makeatletter
\usepackage{amsfonts}
\usepackage{amsmath}

\usepackage{amsfonts}\usepackage{amsthm}\usepackage{color}\usepackage{enumitem}\usepackage{bm}\usepackage{cancel}\usepackage{thmtools}

\usepackage[margin=1.5in]{geometry}

\newcommand{\e}{\varepsilon}

\newcommand{\id}{{\mathrm{id}}}

\newtheorem{theorem}{Theorem}%[section]
\newtheorem{proposition}[theorem]{Proposition}

\newtheorem{claim}[theorem]{Claim}

\theoremstyle{definition}

\newtheorem{remark}[theorem]{Remark}

\theoremstyle{remark}

\begin{document}

\title{An automata group of intermediate growth and exponential activity}

\author{J\'er\'emie Brieussel}
\maketitle
\begin{abstract}
We give a new example of an automata group of intermediate growth. It is generated by an automaton with 4 states on an alphabet with 8 letters. This automata group has exponential activity and its limit space is not simply connected.
\end{abstract}

The growth function $b(n)$ of a group with a finite generating set counts the number of group elements expressible by a word of length at most $n$ in the generators. Up to multiplicative constant in front of the argument, it does not depend upon the choice of generating set. The question of existence of groups with intermediate growth, i.e. growth neither polynomial nor exponential, was asked by Milnor in the 60's and answered positively by Grigorchuk \cite{Grigorchuk83} among some specific automata groups. These groups had already appeared in the work of Aleshin about infinite torsion groups \cite{Aleshin1972}. Grigorchuk's construction has been generalized in a number of ways \cite{Grigorchuk85}, \cite{Grigorchuk1986}, \cite{FabrykowskiGupta1991}, \cite{BartholdiSunik2001}, \cite{bartholdiNonU2003}, \cite{Brieussel2009}, \cite{BartholdiErschler2012}. All these generalizations are groups acting on a rooted tree with bounded activity. 

The activity function $\rm{act}_g(\ell)$ of an automorphism $g$ of a rooted tree counts the number of non-trivial sections at level $\ell$. The activity (polynomial of integer degree $d \geq 0$ or exponential) of an automata group, was defined and classified by Sidki \cite{Sidki2000}. It is the highest growth rate among activity functions of the elements of the group. To the author's knowledge, all groups of intermediate growth constructed by automata are naturally generated by an automaton with bounded, i.e. polynomial of degree 0, activity. We point out that there are many automata of bounded activity generating groups of exponential growth, see for instance \cite[Proposition 7.14]{Bartholdi2017survey}, or of polynomial growth. 

The aim of this note is to exhibit a new automata group with intermediate growth and exponential activity. The automaton generating this group is inspired by a construction due to Wilson of groups of non-uniform exponential growth \cite{Wilson2004}, \cite{Wilson2003}, see Remark~\ref{Wilson}. It should be noted that the activity depends not only on the group but also on the action on the tree. For instance, the Grigorchuk group can be generated by an automaton with exponential activity, see Remark \ref{Grigexp} due to Godin. It is possible that an isomorphic copy of the group described below, or of one of its finite index subgroups, would be generated by an automaton with bounded activity. It seems difficult to rule out this possibility.

Another interesting feature of the group described below, pointed out to the author by Nekrashevych, is that its limit space is not simply connected (in fact it has uncountable fundamental group). The limit space is a renormalized limit of the Schreier graphs of the actions on the levels of the tree, see \cite[Chapter 3]{NekrashevychBook} for a precise definition. All the previously known examples of automata groups of intermediate growth had a dendroid limit space, e.g. an interval for the Grigorchuk group and a dendroid Julia set for the Gupta-Fabrykowsky group, see \cite[Section 6.12]{NekrashevychBook}.

The family of known groups with intermediate growth has been recently drastically extended by Nekrashevych, who constructed the first examples of simple groups of intermediate growth \cite{Nekrashevych2016}. These groups are very different from automata groups since they do not act faithfully on a rooted tree. The present example shows that even among automata groups, in particular among residually finite groups, the zoology of groups with intermediate growth is wider than suggested by the limited list of known examples.

{\bf Notations.} Recall the standard permutational wreath product isomorphism 
\[
\rm{Aut}(T_8) \simeq \rm{Aut}(T_8) \wr_{\{1,\dots,8\}} \mathcal{S}_8
\]
of the group $\rm{Aut}(T_8)$ of automorphisms of an 8-regular rooted tree, where $\mathcal{S}_8$ is the group of permutations of the set $\{1,\dots,8\}$. We write $\langle g_1,\dots,g_8 \rangle \sigma$ for the image of $g$ under this isomorphism and we will identify it with $g$ itself. We call $\sigma$ the permutation of $g$ and $g_i$ the $i$th section of $g$. The reader is refered to \cite{NekrashevychBook} for more about automata groups.

\begin{figure}
	\begin{center}
\begin{tikzpicture}

\draw[thick,blue] (0,0)--(-1.732,0.866);
\draw[blue] (-0.916,0.473)--(-0.716,0.473);
\draw[blue] (-0.916,0.473)--(-0.816,0.3);

\draw[thick,blue] (0,0)--(-1.732,-0.866);
\draw[blue] (-0.916,-0.47)--(-1.136,-0.47);
\draw[blue] (-0.916,-0.47)--(-1.026,-0.643);

\draw[thick,blue] (-1.732,0.866)--(-1.732,-0.866);
\draw[blue] (-1.732,-0.086)--(-1.832,0.086);
\draw[blue] (-1.732,-0.086)--(-1.632,0.086);

\draw[thick,red] (0,0)--(2,0);

\draw[red] (0.6,0)--(0.77,0.1);
\draw[red] (0.6,0)--(0.77,-0.1);
\draw[red] (1.4,0)--(1.23,0.1);
\draw[red] (1.4,0)--(1.23,-0.1);

\draw[thick,blue] (2,0)--(3.732,0.866);
\draw[blue] (2.916,0.47)--(3.136,0.47);
\draw[blue] (2.916,0.47)--(3.026,0.643);

\draw[thick,blue] (2,0)--(3.732,-0.866);
\draw[blue] (2.916,-0.473)--(2.716,-0.473);
\draw[blue] (2.916,-0.473)--(2.816,-0.3);

\draw[thick,blue] (3.732,0.866)--(3.732,-0.866);
\draw[blue] (3.732,0.086)--(3.832,-0.086);
\draw[blue] (3.732,0.086)--(3.632,-0.086);

\draw[thick,red] (3.732,0.866)--(5.732,0.866);
\draw[red] (4.332,0.866)--(4.532,0.966);
\draw[red] (4.332,0.866)--(4.532,0.766);
\draw[red] (5.132,0.866)--(4.962,0.966);
\draw[red] (5.132,0.866)--(4.962,0.766);
\draw[thick,red] (3.732,-0.866)--(5.732,-0.866);

\draw[red] (4.332,-0.866)--(4.532,-0.966);
\draw[red] (4.332,-0.866)--(4.532,-0.766);
\draw[red] (5.132,-0.866)--(4.962,-0.966);
\draw[red] (5.132,-0.866)--(4.962,-0.766);

\draw (-1.732,0.866) node[above]{1};
\draw (-1.732,0.866) node{$\bullet$};
\draw (-1.732,0.866) node[left,red]{$a$};

\draw (-1.732,-0.866) node[below]{2};
\draw (-1.732,-0.866) node{$\bullet$};

\draw (0,0) node[below]{3};
\draw (0,0) node{$\bullet$};

\draw (2,0) node[below]{4};
\draw (2,0) node{$\bullet$};

\draw (3.732,-0.866) node[below]{5};
\draw (3.732,-0.866) node{$\bullet$};

\draw (3.732,0.866) node[above]{6};
\draw (3.732,0.866) node{$\bullet$};

\draw (5.732,0.866) node[above]{7};
\draw (5.732,0.866) node{$\bullet$};
\draw (5.732,0.866) node[right,blue]{$b^{-1}$};

\draw (5.732,-0.866) node[below]{8};
\draw (5.732,-0.866) node{$\bullet$};
\draw (5.732,-0.866) node[right,blue]{$b$};

\end{tikzpicture}
	\end{center}
	\caption{The Schreier graph of the action of the generators $a$ in red and $b$ in blue acting on the first level identified with the 8 letters alphabet.}
	\label{fig1}
\end{figure}
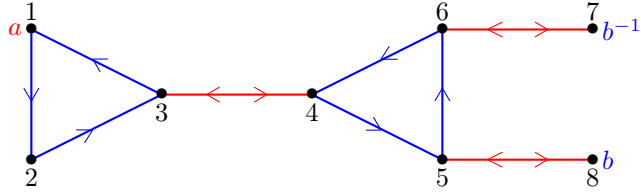

Let $G=\langle a,b,b^{-1}\rangle <\textrm{Aut}(T_8)$ be the automata group over an 8 letters alphabet generated by the three elements recursively defined by
\begin{eqnarray}\label{def}
\begin{array}{lll}a&=&\langle a,\id,\id,\id,\id,\id,\id,\id\rangle(34)(67)(58), \\ b&=&\langle\id,\id,\id,\id,\id,\id,b,b^{-1}\rangle(123)(456),
\\ b^{-1}&=&\langle\id,\id,\id,\id,\id,\id,b^{-1},b\rangle(132)(465).\end{array}
\end{eqnarray}
See Figure \ref{fig1} for an illustration.

\begin{proposition}\label{mainprop}
The automata group $G$ has exponential activity and intermediate growth:
\[
b(n) \le e^{cn^\alpha}, \quad \rm{ where } \quad \alpha=\frac{\ln(8)}{\ln(8)-\ln(\frac{7}{8})} \approx 0.9396.
\]
Moreover, the group $G$ is torsion.
\end{proposition}

The proof of torsion was explained to the author by Laurent Bartholdi who kindly accepted that it is reproduced here. Both proofs of intermediate growth and torsion follow the classical strategy of Grigorchuk \cite{Grigorchuk85}, using self-similarity and strong contraction.

\begin{proof}
The group $G$ obviously has exponential activity because $\textrm{act}_b(\ell)=2^\ell$. First observe that $a^2=b^3=\id$. The key point is that $1$ is a fixed point of the permutation of $a$ and $7,8$ are fixed points of the permutation of $b$. It follows that $G$ is a quotient of the free product $\mathbb{Z}/2\mathbb{Z} \ast \mathbb{Z}/3\mathbb{Z}$. Any minimal length representative word of some element of $G$ has the form
\begin{eqnarray}\label{generic}
w=b^{\varepsilon_1}(ab)^{p_1}(ab^{-1})^{q_1}(ab)^{p_2}(ab^{-1})^{q_2}\dots (ab)^{p_k}(ab^{-1})^{q_k}a^{\e_2},
\end{eqnarray}
with $\e_i \in \{-1,0,1\}$. Moreover, the positive integers $p_i,q_i$ (possibly $p_1$ or $q_k$ are zero), are no more than 7, at the exception of  the following four words: $(ab)^8$, its $a$-conjugate $(ba)^8$ and their inverses $(b^{-1}a)^8$ and $(ab^{-1})^8$. Indeed, using the relation $(ab)^{16}=\id$, we would contradict minimality.

For $i$ from 1 to 8, denote $w_i$ the $\mathbb{Z}/2\mathbb{Z} \ast \mathbb{Z}/3\mathbb{Z}$-free reduction of the $i$th section of $w$. Set $L(w)=\sum_{i=1}^8 |w_i|_a$ where $|w|_a$ is the number of letters $a$ appearing in $w$. Observe that $2|w|_a$ is the length of $w$ up to $\pm1$.

\begin{claim}
\[
L(w) \le \frac{7}{8}|w|_a+1
\]
\end{claim}
To check this claim, observe that each letter $a$ in $w$ contributes at most one letter $a$ in one of the sections $w_i$. Moreover, the equality
\begin{eqnarray}\label{key}
abab^{-1}a=\langle a^2,\id,b,a,b^{-1},b,\id,b^{-1}\rangle(27)(35)(68)
\end{eqnarray}
ensures that whenever the subword $abab^{-1}a$ appears in $w$, the subword $a^2=\id$ appears in a section, reducing its length by two. The claim follows by counting the frequency of non-overlapping such suwords. A worse-case situation is given by $(ab)^7ab^{-1}ab(ab^{-1})^7$. Note that the key point to get (\ref{key}) is that the permutation of $b$ maps 1 to 2, and that both 1 and 2 are fixed by the permutation of $a$.

The claim ensures that the ball $B(n)$ of radius $n$ in the Cayley graph of $G$ embedds into a union of products $B(n_1)\times \dots \times B(n_8)$, where the union is over all possible $n_i$ satisfying $n_1+\dots+n_8 \leq \frac{7}{8}n+c$, of which there is polynomial choice. The precise bound on growth follows by Muchnik-Pak's growth theorem \cite{MuchnikPak2001}.

To prove torsion, proceed by induction on $|w|_a$, were $w$ is a reduced word of the form (\ref{generic}). Assume that $|w|_a>8$. Let $c$ be a cycle of the permutation of $w$, of length $k$. For $i$ a letter in the cycle, $w^k$ fixes $i$ and its $i^{\textrm{th}}$ section is $\prod_{k=0}^{k-1}w_{\sigma^k(i)}$ of length $\leq  \frac{7}{8}|w|_a+1$ by the claim. By induction, some power of $w^k$ is trivial on $c$. As this holds true for each cycle, $w$ is torsion. There remains to check torsion of elements $|w|_a\le 8$, using GAP.

There remains to show the growth of $G$ is not polynomial. But if this were the case, the group would be virtually nilpotent by Gromov's theorem, hence virtually torsion free. Therefore, it is sufficient to prove that $G$ is infinite. This is the case because the section map $g\mapsto g_1$ from $\mathrm{Stab}_G(1)$ to $G$ is onto, where $\mathrm{Stab}_G(1)$ is the stabilizer in $G$ of the vertex $1$ in the tree.  Indeed, its image contains the generators $a$ and $b$, as easily checked using the definition of $a$ and
\[
b^{(abab)}=\langle b,ab^{-1},a,b^{-1}a,b^{-1},\id,b^{-1},ba\rangle (264)(358).
\]
\end{proof}

\begin{remark}\label{Wilson}
The choice of the permutations of $a$ and $b$ follows a construction due to Wilson of groups of non-uniform exponential growth \cite{Wilson2004}, \cite{Wilson2003}. His construction also provides groups of intermediate growth as explained in \cite[Section 6]{Brieussel2009}. One of them is the following automata group $H=\langle a,b'\rangle <\rm{Aut}(T_8)$ generated by $a$ and
\begin{eqnarray}\label{def2}
b'&=&\langle\id,\id,\id,\id,\id,\id,b',\id\rangle(123)(456).
\end{eqnarray}
The group $H$ only has bounded activity. It also has intermediate growth because the relation (\ref{key}) still holds. However, the best known upper bound is $\exp\left(\frac{cn\log\log n}{\log n} \right)$, see \cite[Remark 6.7]{Brieussel2009} which uses an argument of Erschler \cite{Erschler2004boundary}. Perhaps counter-intuitively, augmenting the activity by replacing $b'$ by $b$ permits to obtain a better upper bound on growth, because it gives the new relation $(ab)^{16}=\id$, bounding the values of $p_i,q_i$ in (\ref{generic}).
\end{remark}

\begin{remark}\label{Grigexp}
The first Grigorchuk group of intermediate growth is usually defined via the following automaton of bounded activity: $a=\langle \id,\id\rangle (1 2), b=\langle c,a\rangle, c=\langle d,a\rangle, d=\langle b,\id\rangle$. It was pointed out to the author by Thibault Godin that replacing all the $a$'s in this automaton by the automorphism $a'=\langle a',a'\rangle (12)$, which has exponential activity, yields an isomorphic copy. It is not known if the group of Proposition \ref{mainprop} can be generated by an automaton with bounded activity. 
\end{remark}

{\bf Acknowledgments.} I wish to thank Laurent Bartholdi who explained to me the proof of torsion, Thibault Godin who pointed out Remark \ref{Grigexp}, Slava Grigorchuk for many comments and Volodia Nekrashevych who pointed out the specificity of the limit space. 

\bibliographystyle{alpha}
\bibliography{AutomataT8}

\textsc{\newline J\'er\'emie Brieussel \newline Universit\'e de Montpellier - Institut Montpelli\'erain Alexander Grothendieck }
%\newline Place E. Bataillon cc 051 \newline 34095 Montpellier, France} \newline
\textit{E-mail address:} jeremie.brieussel@umontpellier.fr
\end{document}